\documentclass[reqno]{amsart}

\usepackage[all]{xy}

\numberwithin{equation}{section}

\theoremstyle{plain}
\newtheorem{lemma}{Lemma}[section]
\newtheorem{theorem}[lemma]{Theorem}
\newtheorem{proposition}[lemma]{Proposition}
\newtheorem{corollary}[lemma]{Corollary}

\newtheorem{claim}{Claim}

\theoremstyle{definition}
\newtheorem{definition}[lemma]{Definition}

\newtheorem{example}[lemma]{Example}

\theoremstyle{remark}
\newtheorem{remark}[lemma]{Remark}

\newcommand{\pic}{{\rm Pic\thinspace}}
\newcommand{\bs}{{\rm Bs\thinspace}}
\newcommand{\p}{\mathbb{P}}

\newcommand{\f}{\mathbb{F}}

\newcommand{\oc}{{\mathcal O}}
\newcommand{\ls}{{\mathcal L}}
\newcommand{\ms}{{\mathcal M}}
\newcommand{\h}{\mathcal{H}}
\newcommand{\ns}{{\mathcal N}_b}
\newcommand{\sd}{{\mathcal S}}

\newcommand{\ci}{{\mathcal I}}
\newcommand{\cre}{{\rm Cr\thinspace}}

\newcommand{\rt}{\longrightarrow}
\newcommand{\rmap}{\dashrightarrow}

\begin{document}

\title{On linear systems of $\p^3$ through multiple points}
\author{Cindy De Volder}
\address{
Department of Pure Mathematics and Computeralgebra, Galglaan 2,
\newline B-9000 Ghent, Belgium} \email{cdv@cage.ugent.be}
\thanks{The first author is a Postdoctoral Fellow of
the Fund for Scientific Research-Flanders (Belgium)
(F.W.O.-Vlaanderen)}

\author{Antonio Laface}
\address{
Dipartimento di Matematica, Universit\`a degli Studi di Milano,
Via Saldini 50, \newline 20100 Milano, Italy }
\email{antonio.laface@unimi.it}
\thanks{The second author would like to thank the European Research and
Training Network EAGER for the support provided at Ghent
University. He also acknowledges the support of the MIUR of the
Italian Government in the framework of the National Research
Project ``Geometry in Algebraic Varieties'' (Cofin 2002)}
\begin{abstract}
In this paper we prove a conjecture about the dimension of linear
systems of surfaces of degree $d$ in $\p^3$ through at most eight
multiple points in general position.
\end{abstract}
\maketitle

\section{Introduction}

In this paper we assume the ground field is algebraically closed
of characteristic 0.

The aim of this paper is to evaluate the dimension of linear
systems of surfaces of degree $d$ of $\p^3$ through at most eight
multiple points in general position, i.e.
$\oc_{\p^3}(d)-\sum_{i=1}^8 m_ip_i$. The virtual dimension of the
system is the dimension of the projective space of polynomials of
degree $d$ minus the number of conditions imposed by each point
evaluated independently. It is possible that these conditions are
actually dependent, giving place to the existence of a special
linear system. Recently a conjecture on the structure of special
systems of $\p^3$ has been formulated in~\cite{lu}. In this paper
we provide a proof of this conjecture for these systems. The main
idea is an extension of a procedure introduced in \cite{ha} for
the study of linear systems on $\p^2$ through at most nine
multiple points. By using cubic Cremona transformations of $\p^3$
it is possible to transform a linear system $\ls$ into another one
which is in ``standard form''. The dimensions of the two systems
are the same, while the virtual dimensions may be different and
this difference is measured by proposition~\ref{vir-change}. This
is a completely new phenomenon which does not occur in $\p^2$.
Once the system $\ls$ is in standard form its dimension is related
with that of $\ls_{|Q_i}$, where $Q_i=-\frac{1}{2}K_{X_i}$ is half
of the anti-canonical bundle of the blow-up of $\p^3$ along $i$
points. In this way it is possible to evaluate the speciality of
$\ls$ (theorem~\ref{main}).

The paper is organized as follows: in the first section we recall
some preliminary definitions. In Section $2$ we give a description
of Cremona transformations of $\p^3$ and their action on linear
systems while Section $3$ deals with special linear systems
produced by $(-1)$-curves. In Section $4$ we prove the main
theorem and section $5$ is an appendix on a birational
transformation from $\p^1\times\p^1$ to $\p^2$.

\section{Preliminaries}

We start by fixing some definitions and notations.

Let $Z=\sum m_ip_i$ be a zero-dimensional scheme of general fat
points of $\p^3$; with $\ls=\ls_3(d,\allowbreak
m_1,\allowbreak\ldots,\allowbreak m_r)$ we will denote the linear
system associated to the sheaf $\oc_{\p^3}(d)\otimes\ci_Z$. Given
a linear system $\ls$, by abuse of notation we will denote by
$v(\ls)$ the virtual dimension of the associated sheaf:
\[
v(\ls) := \binom{d+n}{n}-\sum_{i=1}^r\binom{m_i+n-1}{n}-1.
\]
In the same way we will denote by $H^i(\ls)$ the $i$-th cohomology
group of the sheaf associated to $\ls$.

Let $X\stackrel{\pi}\rt\p^3$ be the blow-up of $\p^3$ along
$\{p_1,\ldots,p_r\}$; with abuse of notation we will denote by
$\ls$ the linear system associated to $L=dH-\sum m_iE_i$, where
$H$ is the pull-back of an hyperplane of $\p^3$ and
$E_i=\pi^{-1}(p_i)$. With $\langle h,e_1,\ldots,e_r\rangle$ we
denote a basis of ${\rm\bf A}^2(X)$ where $h$ is the pull-back of
a class of a general line in $\p^3$ and $e_i$ is the class of a
line of $E_i$. The notation
$\ell=\ell_3(\delta,\mu_1,\ldots,\mu_r)$ indicates a curve in
$\p^3$ of degree $\delta$ through $r$ points of multiplicity
$\mu_1,\ldots,\mu_r$ or equivalently a curve $\delta
h-\sum_{i=1}^r\mu_i e_i$ of $X$. The intersection product
$\ell\ls$ is to be intended always on $X$, i.e.
\[
\ell_3(\delta,\mu_1,\ldots,\mu_r)\ls_3(d,m_1,\ldots,m_r)=\delta
d-\sum_{i=1}^r\mu_i m_i.
\]

In what follows we give some definitions in the same spirit of
\cite{ha}.

\begin{definition}
A linear system $\ls=\ls_3(d,m_1,\ldots,m_r)$ is in {\em standard
form} if $2d\geq\sum_{i=1}^4 m_i$ and $m_1\geq\ldots\geq m_r$.
\end{definition}

Let $\sd_i=\ls_3(2,1^i)$ be a linear system of quadrics through
$4\leq i\leq 9$ simple points. We call $\sd_i$ a {\em standard
class}.

\begin{proposition}\label{standard}
A linear system in standard form may be always written in the
following way $\ls = \sd + \sum_{i=4}^a c_i \sd_i$, where the
$c_i$'s are non negative integers and
$\sd=\ls_3(d-2m_4,m_1-m_4,m_2-m_4,m_3-m_4)$.
\end{proposition}
\begin{proof}
Let $a=\max\{i\ |\ m_i>0\}$, if $a\leq 3$ there is nothing to
prove, otherwise consider the system $\ls'=\ls-\sd_a$. Since
$\ls'=\ls_3(d-2,m_1-1,\ldots,m_a-1)$ it follows that also $\ls'$
is in standard form. Proceeding in the same way on $\ls'$, after a
finite number of steps one obtains a linear system $\sd$ through
at most three points. In this way one obtains that $c_a = m_a$ and
$c_{i}=m_i-m_{i+1}$ for $i < a$ while $\ls^{(i)}=\ls-\sum_{j=i}^a
c_j\sd_j$ is given by $\ls_3(d-2m_i,m_1-m_i,\ldots,m_{i+1}-m_i)$.
\end{proof}

Let $Q$ be a quadric, with $\ls_Q(a,b,m_1,\ldots,m_r)$ we will
mean the system $|\oc_Q(a,b)|$ through $r$ general points of
multiplicities $m_1,\ldots,m_r$.

\section{Cubic Cremona transformations of $\p^3$}

In this section we focus our attention on a class of cubic Cremona
transformations of $\p^3$. Consider the system $\ls_3(3,2^4)$, by
putting the four double points in the fundamental ones, the
associated rational map is given by:
\begin{equation}\label{cubic}
\cre: (x_0:x_1:x_2:x_3) \rmap
(x_0^{-1}:x_1^{-1}:x_2^{-1}:x_3^{-1}).
\end{equation}
This birational map induces an action on the picard group of the
blow-up $X$ of $\p^3$ along the four points which can be described
in the following way:
\begin{proposition}[\cite{lu}]\label{cre-surfaces}
The action of transformation~\ref{cubic} on
$\ls=\allowbreak\ls_3(d,\allowbreak m_1,\allowbreak \ldots,m_r)$
is given by:
\begin{eqnarray}\label{cre-a1}
\cre(\ls) & := & \ls_3(d+k,m_1+k,\ldots, m_4+k,m_5,\ldots,m_r),
\end{eqnarray}
where $k=2d-\sum_{i=1}^4m_i$.
\end{proposition}

Observe that $\dim \cre(\ls) = \dim \ls$ but in general the
virtual dimensions of the two systems may be different.

\begin{proposition}[\cite{lu}]\label{vir-change}
Let $\ls=\ls_3(d,m_1,\ldots,m_r)$ be a linear system such that $2d
\geq m_i + m_j + m_k$ for any choice of
$\{i,j,k\}\subset\{1,2,3,4\}$ then
\begin{equation}\label{vc}
v(\cre(\ls)) - v(\ls) = \sum_{t_{ij} \geq
2}\binom{1+t_{ij}}{3}-\sum_{t_{ij}\leq -2}\binom{1-t_{ij}}{3},
\end{equation}
where $t_{ij}=m_i+m_j-d$.
\end{proposition}

\begin{corollary}
Under the same assumptions of Proposition~\ref{vir-change}, if the
degree of $\cre(\ls)$ is smaller than that of $\ls$, then
$v(\cre(\ls))\geq v(\ls)$.
\end{corollary}
\begin{proof}
The difference between the degree of $\cre(\ls)$ and that of $\ls$
is equal to $k = 2d - \sum_{i=1}^4m_i$. From $2d <
\sum_{i=1}^4m_i$ we deduce that, if $t_{12}\geq 2$ then $d-m_3-m_4
< m_1+m_2-d$ which is equivalent to $-t_{34} < t_{12}$. The same
holds for each $t_{ij}$ such that $t_{ij}\geq 2$, hence the right
side of equation~\ref{vc} is non negative.
\end{proof}

%The birational map~\ref{cubic} gives an action also on ${\rm\bf
%A}^2(X)$ which can be described as follows:

%\begin{proposition}[\cite{lu}]\label{cre-curves}
%The action of transformation~\ref{cubic} on
%$\ell=\allowbreak\ell_3(\delta,\allowbreak \mu_1,\allowbreak
%\ldots,\mu_r)$ is given by:
%\begin{eqnarray}\label{cre-a2}
%\cre(\ell) & := & \ell(\delta+2h,\mu_1+h,\ldots,
%\mu_4+h,\mu_5,\ldots,\mu_r),
%\end{eqnarray}
%where $h=\delta-\sum_{i=1}^4\mu_i$. \\
%\end{proposition}

%From proposition~\ref{cre-surfaces} and~\ref{cre-curves} it
%follows immediately that
%$
%\cre(\ls)\cre(\ell)=\ls\ell.
%$

\section{$(-1)$-curves and special systems}

Starting from the results of the preceding section it is possible
to define a class of special linear systems.

\begin{definition}
A curve $C\in\ell$ is called $(-1)$-curve if $\ell$ is obtained by
applying a finite set of Cremona transformations on the system
$\ell_3(1,1^2)$.
\end{definition}

\begin{example}
Given six general points of $\p^3$ there exists a unique rational
normal curve through them. This curve is an element of
$\ell_3(3,1^6)$ which may be obtained as
$\cre(\ell_3(1,0^4,1^2))$, i.e. as the Cremona transformation of a
line through two points.
\end{example}

\begin{proposition}\label{speciality}
Let $\ls=\ls_3(d,m_1,\ldots,m_r)$ be a linear system and
$\ell_1,\ldots,\ell_n$ be a set of $(-1)$-curves such that
$\ell_i\ls=-t_i\leq -2$ for $i=1\ldots,n$. Then
\[
\dim\ls - v(\ls) \geq \sum_{i=1}^n\binom{t_i+1}{3} -
h^2(\ls\otimes\ci_{\ell_1}^{t_1}\ldots\otimes\ci_{\ell_n}^{t_n}),
\]
where with $\ls\otimes\ci_{\ell_i}^{t_i}$ we mean the sheaf
$\oc_{\p^3}(d)\otimes\ci_Z\otimes\ci_{\ell_i}^{t_i}$.
\end{proposition}
\begin{proof}[Idea of the proof]
For a complete proof of this proposition see~\cite{lu}. \\
By definition, each one of the $\ell_i$ is given by a smooth curve
$C_i\subset X$ such that
$N_{C_i|X}\cong\oc_{\p^1}(-1)\oplus\oc_{\p^1}(-1)$. Consider the
blow-up $Y\stackrel{p}\rt X$ along the curves $C_1,\ldots,C_n$,
the exceptional divisors are quadrics $Q_1,\ldots,Q_n$. From the
evaluation of the tautological line bundle associated to the
blow-up of $C_i$ and from the intersection $C_i\ls=-t_i$ one
obtains:
\begin{eqnarray*}
{Q_i}_{|Q_i} & \cong & \oc(-1,-1) \\
p^*\ls_{|Q_i} & \cong & \oc(-t_i,0).
\end{eqnarray*}
these formulas imply that $t_iQ_i\subseteq\bs (p^*\ls)$ and that
\[
\chi(p^*\ls)=\chi(p^*\ls-t_iQ_i)+\sum_{k=0}^{t_i}\chi(\oc(k-t_i,k)).
\]
An easy calculation shows that the last sum is equal to
$\binom{t_i+1}{3}$, hence applying this procedure to each one of
the $C_i$ one obtains:
\[
\chi(p^*\ls)=\chi(p^*\ls-\sum_{i=1}^nt_iQ_i)+\sum_{i=1}^{n}\binom{t_i+1}{3}.
\]
This, together with the fact that
$h^0(p^*\ls)=h^0(p^*\ls-\sum_{i=1}^nt_iQ_i)$ proves the thesis.
\end{proof}

\begin{example}\label{h2non0}
The system $\ls=\ls_3(3,3^3)$ has $v(\ls)=-11$ while $\dim\ls =0$
since it consists of three times the plane $\ls_3(1,1^3)$. For
each line $\ell_{i}$ through two of the three points, we have
$\ls\cdot\ell_{i}=-3$, hence the speciality is greater then or
equal to
$3\binom{4}{3}-h^2(\ls\otimes\ci_{\ell_{1}}^3\otimes\ci_{\ell_{2}}^3\otimes\ci_{\ell_{3}}^3)$.
So the $h^2$ is equal to $1$.
\end{example}

The preceding proposition allow us to give an estimate of the
speciality of a given linear system $\ls$. In particular consider
a system in standard form, then we have the following:

\begin{proposition}\label{h2}
Let $\ls=\ls_3(d,m_1,\ldots,m_8)$ be a linear system in standard
form, let $t_i := m_1+m_i-d$ for $i\geq 2$ and $t_1 := m_2+m_3-d$,
then the following holds:
\[
h^1(\ls)\geq \sum_{t_i\geq 2}\binom{t_i+1}{3}.
\]
\end{proposition}
\begin{proof}
Let $l_i$ be a line through two multiple points of $\ls$, we say
that the $1$-dimensional scheme $\Gamma=\sum t_il_i$ is {\em
associated} to $\ls$ if $l_i\ls=-t_i<0$. By
proposition~\ref{speciality} it is sufficient to prove that
$h^2(\ls\otimes\ci_{\Gamma})=0$. Consider the system given by the
plane through the first three points $\h=\ls_3(1,1^3)$ and let
$\Gamma=\Gamma'+\Gamma''$ where $\Gamma'=l_2+l_3$, then one has
the exact sequence
\begin{equation}\label{plane}
\xymatrix@1{ 0 \ar[r]  & (\ls-\h)\otimes\ci_{\Gamma''} \ar[r] &
\ls\otimes\ci_{\Gamma} \ar[r] & (\ls\otimes\ci_{\Gamma})_{|W}
\ar[r] & 0, }
\end{equation}
where $W\in\h\otimes\ci_{\Gamma'}$ is the plane through the first
three points.
\begin{claim}\label{h2-sequence}
In the preceding sequence $h^2((\ls\otimes\ci_{\Gamma})_{|W})=0$
and $\ls-\h$ is a system in standard form unless
$\ls=\ls_3(2m+t,m+2t,m^5,m_7,m_8)$, with $1\leq t\leq m$.
\end{claim}
Observe that $\Gamma''$ is associated to $\ls-\h$, since this
system is equal to  $\ls_3(d-1,m_1-1,m_2-1,m_3,\ldots,m_8)$. So,
if $\ls-\h$ is still in standard form, then after sorting the new
multiplicities we can proceed by induction on the length of
$\Gamma$; once we obtain a system $\bar{\ls}$ for which the
associated $\bar{\Gamma}$ is empty, we know that
$h^2(\bar{\ls})=0$. \\
Suppose now that $\ls-\h$ is not in standard form, then the system
must be the one of the claim. The worst possibility is that
$m_7=m_8=m$, so we fix our attention on the system
$\ls_3(2m+t,m+2t,m^7)$. By abuse of notation, we still call this
system $\ls$. Applying sequence~\ref{plane} one time we obtain the
system $\ls-\h=\ls_3(2m+t-1,m+2t-1,(m-1)^2,m^5)$. Now, let
$\h_1=\ls_3(2,2,1^5)$ be the quadric cone through the point of
multiplicity $m+2t-1$ and the five points of multiplicity $m$ and
let $W_1\in\h_1$, $\ls_1=\ls-\h$, and
$\Gamma_1=\Gamma_1'+\Gamma_1''$ with $\Gamma_1'=\sum_{i=4}^8l_i$
and $\Gamma_1''=(t-1)\sum_{i=2}^8l_i$. Observe that $\Gamma_1$ is
associated to $\ls_1$, while $\Gamma_1''$ is associated to
$\ls_1-\h_1$. The last system is equal to
$\ls_3(2m'+t',m'+2t',m'^7)$ where $m'=m-1$ and $t'=t-1$, hence by
using sequence~\ref{plane} with $\h_1$ instead of $\h$ and by the
following
\begin{claim}\label{cone}
With the preceding notation
$h^2((\ls_1\otimes\ci_{\Gamma'_1})_{|W_1})=0$,
\end{claim}
we can make induction on $t$ and proving that
$h^2(\ls\otimes\ci_{\Gamma})=0$.
\end{proof}

\begin{proof}[Proof of Claim~\ref{h2-sequence}]
First of all observe that
\[
(\ls\otimes\ci_{\Gamma})_{|W} =
\ls_2(d-t_2-t_3,m_1-t_2-t_3,m_2-t_2,m_3-t_3).
\]
Since $d-t_2-t_3 = (d-m_1)+(2d-m_1-m_2-m_3)\geq 0$ then the $H^2$
of this linear system vanishes. The sequence of multiplicities of
$\ls-\h$ is
\[
m_1-1,m_2-1,m_3-1,m_4,\ldots,m_8
\]
and the degree is $d-1$. This means that the system is still
stable if at least $m_1-1$ and $m_2-1$ are between the biggest
four. If this is not the case then $m_2-1\leq m_6$ and this
implies that $m_2=\cdots =m_6$, call this number $m$. Since
$\ls-\h$ is not in standard form we have that $2(d-1)< m_1-1+3m$,
but we also know that $2d\geq m_1+3m$ and this implies that
$2d=m_1+3m$. Let $t=t_2=\cdots =t_6$, then the preceding equation
gives $m_1=m+2t$ and $d=2m+t$.

\end{proof}

\begin{proof}[Proof of Claim~\ref{cone}]
The blow-up of a quadric cone along the vertex is an $\f_2$
surface, hence the strict transform $\tilde{W}_1$ of $W_1$ of the
blow-up of $\p^3$ along the six points is a blow-up of an $\f_2$
surface along five points. The vanishing of
$h^2((\ls_1\otimes\ci_{\Gamma'_1})_{|W_1})$ is equivalent to the
vanishing of
$h^2((\tilde{\ls}_1\otimes\ci_{\tilde{\Gamma}_1})_{|\tilde{W}_1})$
where $\tilde{\ls}_1$ is the strict transform of $\ls_1$ and
$\tilde{\Gamma}_1$ is a $1$-dimensional subscheme of $\tilde{W}_1$
corresponding to the strict transform of $\Gamma'_1$. A basis for
the picard group of ${W}_1$ may be written as $\langle
f,c,e_1,\ldots,e_5\rangle$ where $f^2=0,\ fc=1,\ c^2=-2$ and the
$e_i$ are $(-1)$-curves of the blow-up. With this notation, the
hyperplane section of $\tilde{W}_1$ is given by $c+2f$ (it is very
ample outside $c$ which is contracted to the vertex of the cone).
Instead of the system $\tilde{\ls}_1$ we can consider
$|(2m+t-1)H-(m+2t-1)E-\sum_{i=1}^5mE_i|$, since the two
exceptional divisors of multiplicity $(m-1)$ have no intersection
with $\tilde{W}_1$. In this way we have:
\begin{eqnarray*}
(\tilde{\ls}_1\otimes\ci_{\tilde{\Gamma}_1})_{|\tilde{W}_1} & = &
(2m+t-1)(c+2f)-(m+2t-1)c \\
 & & -\sum_{i=1}^5me_i-\sum_{i=1}^5t(f-e_i) \\
& = & (m-t)c + (4m-3t-2)f -\sum_{i=1}^5(m-1)e_i.
\end{eqnarray*}
Since $t\leq m$ this implies that both the coefficients of $c$ and
$f$ are non-negative and this, by adjunction, implies the
vanishing of the $H^2$ of this divisor.
\end{proof}

\section{Linear systems through at most 8 points}

In what follows we will denote by $\ls$ a linear system of type
$\ls_3(d,\allowbreak m_1,\allowbreak \ldots,\allowbreak m_8)$ in
standard form and let $\sd+\sum_{i=4}^a c_i\sd_i$ be its
decomposition.

\begin{lemma}\label{three}
A linear system $\ls=\ls_3(d,m_1,m_2,m_3)$ is empty if and only if
and only if there exists $i\in\{1,2,3\}$ such that $d < m_i$.
\end{lemma}
\begin{proof}
One part of the proof is trivial, since there are no surfaces of
degree $d$ with a singularity of multiplicity greater than $d$. On
the other hand, if all the $m_i$ are equal to $d$, then the system
$\ls_3(d,d^3)$ is non-empty since it contains $d$ times the plane
$\ls_3(1,1^3)$.
\end{proof}

\begin{lemma}\label{restricted}
If $d\geq m_1$ then $h^1(\ls_{|Q_a})=0$ and $h^0(\ls_{|Q_a})>0$
where $Q_a\in\sd_a$.
\end{lemma}
\begin{proof}
The system on the quadric
$\ls_{|Q_a}=\ls_{Q_a}(d,d,m_1,\ldots,m_a)$ is equivalent (by
\ref{app}) to $\ls_2(2d-m_1,(d-m_1)^2,m_2,\ldots,m_a)$. This is a
plane system through at most $9$ points. First of all we want to
see if it is possible to put it in standard form (i.e. a plane
Cremona transformation can not decrease its degree). Observe that
the three bigger multiplicities may be: $\{d-m_1,d-m_1,m_2\}$,
$\{m_2,m_3,d-m_1\}$, $\{m_2,m_3,m_4\}$. In the first and third
case it is obvious that $2d-m_1$ is greater then or equal to the
sum of these multiplicities (since the system $\ls$ is in standard
form). In the second case the inequality $2d-m_1\geq
d-m_1+m_2+m_3$ is true only if $d\geq m_2+m_3$, so we may assume
that $d=m_2+m_3-t$ where $t\geq 1$. After applying a Cremona
transformation to this system we obtain the following
$\ls_2(d',m_1',m_2',m_3',d-m_1,m_4,\ldots,m_a)$ where
$d'=3d-m_1-m_2-m_3,\ m_1'=d-m_3,\ m_2'=d-m_2,\
m_3'=2d-m_1-m_2-m_3$. The bigger multiplicities of this system are
$\{d-m_3,\ d-m_2,\ d-m_1\}$ since $d-m_1\geq m_4$ by assumption
and $2d-m_1-m_2-m_3=d-m_1-t<d-m_1$. This implies that the system
is in standard form ($d'=d-m_3+d-m_2+d-m_1$). \\
So we proved that after a quadratic transformation of $\p^2$, the
system $\ls_2(2d-m_1,(d-m_1)^2,m_2,\ldots,m_a)$ becomes $\ms$
which is in standard form. By~\cite{ha} this implies that
$h^0(\ms) > 0$. The intersection $\ms .K_{\tilde{\p}^2} = -4d +
m_1\ldots + m_a$, where $\tilde{\p}^2$ is the blow-up of $\p^2$
along the $a+1$ points, is non-positive, so  $h^1(\ms)=0$
(by~\cite{ha}).
\end{proof}

\begin{theorem}\label{main}
A linear system $\ls(d,m_1,\ldots,m_a)$ (with $a\leq 8$) in
standard form is special if and only if $d\leq m_1+m_2 -2$ and its
dimension is given by
\[
\dim\ls = v(\ls) + \sum_{t_{i}\geq 2}\binom{t_i+1}{3},
\]
where $t_1:=m_2+m_3-d$ and $t_i:=m_1+m_i-d$ for $i\geq 2$.
\end{theorem}
\begin{proof}
Since $\ls$ is in standard form, by \ref{standard} it can be
written as $\ls = \sd + \sum_{i=4}^a c_i \sd_i$. We will
distinguish two cases:

If $\sd\neq\emptyset$ it follows immediately that $h^0(\ls)>0$
since the system may be written as a sum of effective ones. In
order to see that $h^1(\ls)=h^1(\sd)$, consider the exact
sequence:

\begin{equation}\label{seq}
\xymatrix@1{ 0 \ar[r]  & \ls'-\sd_i \ar[r] & \ls' \ar[r]  &
\ls'_{|Q_i} \ar[r] &  0, }
\end{equation}

where $\ls'$ is obtained from $\ls$ by subtracting some of the
$\sd_j$. The degree of $\ls'$ is greater than or equal to its
first multiplicity, otherwise it would be empty but this is not
possible by the assumption on $\sd$. By lemma \ref{restricted}
$h^1(\ls'_{|Q_i})=0$, which implies that $h^1(\ls')\leq
h^1(\ls'-\sd_i)$. This gives the following
\[
h^1(\ls)\leq h^1(\ls-\sd_a)\leq\ldots\leq h^1(\sd).
\]
The speciality of $\sd$ may be given only by the lines $\langle
p_i, p_j\rangle$ with $1\leq i<j\leq 3$ if $\delta \leq
\mu_i+\mu_j -2$. From the equality $d-m_i-m_j=\delta-\mu_i-\mu_j$
one has the same speciality also for $\ls$ and this implies that
$h^1(\ls)\geq h^1(\sd)$.

From proposition~\ref{three} we know that $\sd=\emptyset$ if and
only if the degree of the system is less then one of its
multiplicities. By proposition~\ref{standard} this means $d-2m_4 <
m_1-m_4$. Recall that $t_4=m_1+m_4-d$, since $\ls$ is in standard
form the inequality $2d=d+m_1+m_4-t_4\geq m_1+m_2+m_3+m_4$ gives
$d\geq m_2+m_3+t_4$. This implies that $d\geq m_i+m_j$ for each
$i\geq 2,\ j\geq 3,\ i\neq j$. So the speciality of $\ls$ coming
from lines is due only to $\langle p_1, p_i\rangle$ with $i\geq
2$. Recall that $t_i=m_1+m_j-d$ and let $b=\max\{i\ |\ t_i\geq
1\}$. Observe that by definition $t_i-t_j=m_i-m_j$ and that $b\geq
4$. By proposition~\ref{h2} we have the following inequality:
\begin{equation}\label{ineq}
h^1(\ls) \geq \sum_{i=2}^b \binom{t_i+1}{3}.
\end{equation}

Consider the system
\[
\ns = \sd+\sum_{i=4}^{b-1}c_1\sd_i+(t_b-1)\sd_b,
\]
We now use the following:
\begin{claim}\label{c1}
Under these assumptions, $h^0(\ns)=0$ and $h^1(\ns)=\sum_{i=2}^b
\binom{t_i+1}{3}$. Furthermore, if $d''$ and $m_1''$ are,
respectively, the degree and the first multiplicity of the system
$\ns+\sd_b$ then $d''=m_1''$.
\end{claim}

Consider the exact sequence
\[
\xymatrix@1{ 0 \ar[r]  & \ns \ar[r] & \ns+\sd_b \ar[r]  &
{\ns+\sd_b}_{|Q_b} \ar[r] &  0, }
\]
by claim \ref{c1} and lemma \ref{restricted} we know that
$h^1({\ns+\sd_b}_{|Q_b})=0$, this implies that

\[
h^1(\ls)\leq h^1(\ls-\sd_a)\leq\ldots \leq h^1(\ns).
\]
This together with lemma \ref{restricted} and inequality
\ref{ineq} implies that $h^1(\ls) = \sum_{i=2}^b
\binom{t_i+1}{3}$.

\end{proof}

\begin{proof}[Proof of Claim \ref{c1}]
The system $\ns$ is given by
$\ls_3(d-2m_b+2t_b-2,m_1-m_b+t_b-1,\ldots,m_{b-1}-m_b+t_b-1,t_b-1)$.
So $d'=d-2m_b+2t_b-2=m_1-m_b+t_b-2=m_1'-1$ implies that
$h^0(\ns)=0$ and $h^1(\ns)=-\chi(\ns)$. The last quantity is
\[
\chi(\ns) = \binom{d'+3}{3} - \binom{m_1'+2}{3}-\sum_{i=2}^b
\binom{m_i'+2}{3}.
\]
Since $d'=m_1'-1$ the first two terms vanish and from
$m_i'=m_i-m_b+t_b-1=t_i-1$ one obtains $\chi(\ns)=\sum_{i=2}^b
\binom{t_i+1}{3}$ which proves the first part of the claim. For
the second part observe that $d''=d'+2$ and $m_1''=m_1'+1$.

\end{proof}

All these results may be summarized in the following procedure
which allows us to evaluate the dimension of a linear system and
its speciality.

\begin{remark}
Take a linear system $\ls=\ls_3(d,m_1,\ldots,m_8)$ and let
$v:=v(\ls)$.
\begin{itemize}
\item[1 -] Sort the multiplicities in descending order.
\item[2 -] If $2d-m_1-m_2-m_3<0$ remove the plane $\h$ through the
first three points, redefine $\ls$ as $\ls-\h$ and goto step $1$.
\item[3 -] If $2d-m_1-m_2-m_3-m_4<0$ make a cubic Cremona
transformation, redefine $\ls$ as $\cre(\ls)$ and goto step $1$.
\item[4 -] Evaluate $d=\dim \ls$ with theorem~\ref{main}.
\end{itemize}
\end{remark}

\section{Appendix on a birational map $\p^1\times\p^1 \rmap \p^2$}

In this section we consider a birational map $\varphi:
\p^1\times\p^1 \rmap \p^2$ given by blowing up a point on the
quadric and contracting the strict transforms of the two lines
through it. In this way it is possible to give a correspondence
between linear systems through fat points on the quadric and those
on the projective plane. Let us consider a linear system
$\ls_Q(a,b,m)$, i.e. a system of curves on the quadric $Q$ of kind
$\oc(a,b)$ through one point $p$ of multiplicity $m$. Blowing up
the quadric at $p$, the strict transform of the preceding system
is given by $af_1+bf_2-me_p$, where $f_1,f_2$ are the pull-back of
the two rulings of $Q$ and $e_p$ is the exceptional divisor of the
blow-up $\tilde{Q}$. A base change in $\pic(\tilde{Q})$ allows us
to write this divisor as
$(a+b-m)(f_1+f_2-e_p)-(b-m)(f_1-e_p)-(a-m)(f_2-e_p)$. Since the
divisors $f_i-e_p$ are $(-1)$-curves, they can be contracted
giving a linear system on $\p^2$ of degree $a+b-m$ through two
points of multiplicity $b-m$ and $a-m$. This implies that the map
$\varphi$ induces the following correspondence:
\begin{equation}\label{app}
\ls_Q(a,b,m,m_1,\ldots,m_r) \longleftrightarrow
\ls_2(a+b-m,b-m,a-m,m_1,\ldots,m_r).
\end{equation}
It is an easy computation to verify that the virtual dimensions of
the two systems are the same.

\bibliographystyle{plain}

\begin{thebibliography}{May72}

\bibitem{ha}
Brian Harbourne.
\newblock The geometry of rational surfaces and {H}ilbert functions of points
  in the plane.
\newblock In {\em Proceedings of the 1984 Vancouver conference in algebraic
  geometry}, volume~6 of {\em CMS Conf. Proc.}, pages 95--111, Providence, RI,
  1986. Amer. Math. Soc.

\bibitem{lu}
Antonio Laface and Luca Ugaglia.
\newblock On a class of special linear systems of $\p^3$.
\newblock {\em math.AG/03}, preprint.

\end{thebibliography}

\end{document}